\documentclass[11pt]{article}

\usepackage{amscd,amsmath, amssymb, fancyhdr}

\newcommand{\version}{version 3.0,\ \   26.03.2008}

\numberwithin{equation}{section}

\newcommand{\CC}{\mathbb{C}}

\newcommand{\RR}{\mathbb{R}}
\newcommand{\ZZ}{\mathbb{Z}}

\newcommand{\al}{\alpha}
\newcommand{\be}{\beta}

\newcommand{\f}{\varphi}

\newcommand{\abs}[1]{\vert #1\vert}

\def\eqref#1{(\ref{#1})}

\newcommand{\arrow}{{\:\longrightarrow\:}}
\newcommand{\Z}{{\Bbb Z}}
\newcommand{\C}{{\Bbb C}}
\newcommand{\R}{{\Bbb R}}

\newcommand{\6}{\partial}
\def\1{\sqrt{-1}\:}

\renewcommand{\tilde}{\widetilde}
\renewcommand{\bar}{\overline}
\renewcommand{\phi}{\varphi}
\renewcommand{\epsilon}{\varepsilon}
\renewcommand{\geq}{\geqslant}

\newcommand{\im}{\operatorname{im}}
\newcommand{\Lie}{\operatorname{Lie}}
\newcommand{\End}{\operatorname{End}}

\newcommand{\Aut}{\operatorname{Aut}}

\newcommand{\const}{\operatorname{const}}

\newcounter{Mycounter}[section]
\newcounter{lemma}[section]
\renewcommand{\thelemma}{{Lemma \thesection.\arabic{lemma}}}
\newcommand{\lemma}{%
     \setcounter{lemma}{\value{Mycounter}}
     \refstepcounter{lemma}
     \stepcounter{Mycounter}
     {\noindent \bf \thelemma.\ }}

\newcounter{claim}[section]
\newcounter{sublemma}[section]
\newcounter{corollary}[section]
\renewcommand{\thecorollary}{{Corollary
\thesection.\arabic{corollary}}}
\newcommand{\corollary}{%
     \setcounter{corollary}{\value{Mycounter}}
     \refstepcounter{corollary}
     \stepcounter{Mycounter}
     {\noindent \bf \thecorollary.\ }}

\newcounter{theorem}[section]
\renewcommand{\thetheorem}{{Theorem \thesection.\arabic{theorem}}}
\newcommand{\theorem}{%
     \setcounter{theorem}{\value{Mycounter}}
     \refstepcounter{theorem}
     \stepcounter{Mycounter}
     {\noindent \bf \thetheorem.\ }}

\newcounter{conjecture}[section]
\renewcommand{\theconjecture}{{Conjecture
\thesection.\arabic{conjecture}}}
\newcommand{\conjecture}{%
     \setcounter{conjecture}{\value{Mycounter}}
     \refstepcounter{conjecture}
     \stepcounter{Mycounter}
     {\noindent \bf \theconjecture.\ }}

\newcounter{proposition}[section]
\renewcommand{\theproposition}
       {{Proposition \thesection.\arabic{proposition}}}
\newcommand{\proposition}{%
     \setcounter{proposition}{\value{Mycounter}}
     \refstepcounter{proposition}
     \stepcounter{Mycounter}
     {\noindent \bf \theproposition.\ }}

\newcounter{definition}[section]
\renewcommand{\thedefinition}
       {{Definition~\thesection.\arabic{definition}}}
\newcommand{\definition}{%
     \setcounter{definition}{\value{Mycounter}}
     \refstepcounter{definition}
     \stepcounter{Mycounter}
     {\noindent \bf \thedefinition.\ }}

\newcounter{example}[section]
\newcounter{remark}[section]
\renewcommand{\theremark}{{Remark \thesection.\arabic{remark}}}
\newcommand{\remark}{%
     \setcounter{remark}{\value{Mycounter}}
     \refstepcounter{remark}
     \stepcounter{Mycounter}
     {\noindent \bf \theremark.\ }}

\newcounter{problem}[section]
\renewcommand{\theproblem}{{Problem \thesection.\arabic{problem}}}
\newcommand{\problem}{%
     \setcounter{problem}{\value{Mycounter}}
     \refstepcounter{problem}
     \stepcounter{Mycounter}
     {\noindent \bf \theproblem.\ }}

\newcounter{question}[section]
\renewcommand{\leftmark}%
{{\scriptsize  Morse-Novikov cohomology of LCK manifolds}}

\makeatletter

\@addtoreset{equation}{section}
\@addtoreset{footnote}{section}
\makeatother

\def\blacksquare{\hbox{\vrule width 5pt height 5pt depth 0pt}}
\def\endproof{\blacksquare}

\begin{document}
\begin{center}
{\LARGE\bf Morse-Novikov
cohomology of locally conformally K\"ahler manifolds\\[3mm]
}

Liviu Ornea\footnote{Liviu Ornea is partially supported by grant
2-CEx-06-11-22/25.07.2006.} and Misha Verbitsky.

\end{center}

{\small
\hspace{0.15\linewidth}
\begin{minipage}[t]{0.7\linewidth}
{\bf Abstract} \\
A locally conformally K\"ahler (LCK) manifold is a
complex manifold admitting a K\"ahler covering, with
the monodromy acting on this covering by holomorphic homotheties.
We define three cohomology invariants, the Lee class,
the Morse-Novikov class, and the Bott-Chern class,
of an LCK-structure. These invariants play together the
same role as the K\"ahler class in K\"ahler geometry.
If these classes coincide for two LCK-structures,
the difference between these structures can be
expressed by a smooth potential, similar to the
K\"ahler case. We show that the Morse-Novikov class
and the Bott-Chern class of a Vaisman manifold vanish.
Moreover, for any LCK-structure on a manifold,
admitting a Vaisman structure, we prove
that its Morse-Novikov class vanishes. We show
that a compact LCK-manifold $M$ with vanishing Bott-Chern
class admits a holomorphic embedding into a Hopf manifold,
if  $\dim_\C M \geq 3$, a result which parallels the
Kodaira embedding theorem.
\end{minipage}
}

\tableofcontents


\section{Introduction}


\subsection{Cohomological invariants of LCK-structures}

A locally conformally K\"ahler (LCK) manifold is a complex
manifold $M$ (in this paper we shall usually assume $\dim_\C M\geq 3$),
equipped with a Hermitian metric
$\omega$ in such a way that a certain covering $\tilde M$
of $M$ is K\"ahler, and its K\"ahler metric is conformal
to the pullback of $\omega$ by the covering map. In this case, we have
$d\omega=\omega\wedge\theta$, and the 1-form $\theta$
is called {\bf the Lee form} of $M$. It is easy to see
that $\theta$ is closed: indeed,
\[
0=d^2\omega = d\theta \wedge \omega +
\theta\wedge \theta \wedge \omega = d\theta \wedge \omega,
\]
but the multiplication by $\omega$ induces an embedding
$\Lambda^2(M) \arrow \Lambda^4(M)$ if
$\dim_\C M \geq 3$.\footnote{When $\dim_\C M =2$,
one defines LCK-manifolds in the same way, but
$d\theta=0$ should be assumed as a part of
the definition.} The LCK-manifolds are
often considered up to conformal equivalence;
indeed, a manifold which is conformally
equivalent to an LCK-manifold is itself
locally conformally K\"ahler.

If one performs a conformal change, $\omega_1 =
e^f\omega$, the Lee form $\theta$ changes to
$\theta_1 = \theta + df$. The cohomology class
$[\theta]\in H^1(M)$ is an important invariant
of an LCK-manifold. Clearly, $[\theta]\in H^1(M)$
vanishes if and only if $\omega$ is
conformally equivalent to a K\"ahler structure.
In this case $(M, \omega)$ is called
{\bf globally conformally K\"ahler}.

Another, more subtle, invariant
of an LCK-manifold is called {\bf the Mor\-se\--Novikov
class } of $[\omega]\in H^2_\theta(M)$, defined
as follows. Recall that the Morse--Novikov cohomology,
also known as Lichnerowicz cohomology
(defined independently by Novikov and Lichnerowicz
in \cite{_Lichnerowicz_} and \cite{Novikov})
is the cohomology of the complex
\begin{equation}\label{_MN_Equation_}
\Lambda^0(M) \stackrel{d-\theta}\arrow \Lambda^1(M) \stackrel{d-\theta}\arrow
\Lambda^2(M)\stackrel{d-\theta}\arrow \cdots
\end{equation}
with $\theta$ a closed 1-form (see Section \ref{_MN_Section_}
for a more detailed exposition). It is well--known that
the cohomology $H^*_\theta(M)$ of \eqref{_MN_Equation_}
is naturally identified with the cohomology of the local system
given by the character $\pi_1(M)\arrow \R^{>0}$
associated with the cohomology class $[\theta]\in H^1(M)$.

An LCK-form $\omega$ on an LCK-manifold satisfies
$d\omega=\omega\wedge\theta$, therefore it is
$(d-\theta)$-closed. The cohomology class
$[\omega]\in H^2_\theta(M)$ is called
{\bf the Morse--Novikov class of the LCK-manifold}.
It is an invariant of the LCK-manifold,
roughly analogous to the K\"ahler
class on a K\"ahler manifold.

The third, even more subtle, cohomology invariant
of an LCK-manifold has a complex-analytic nature;
it is a Morse--Novikov version of the usual Bott--Chern
class of a closed $(1,1)$-form.

On a compact K\"ahler manifold, an exact $(p,q)$-form
$\eta$ satisfies
\begin{equation}\label{_harm_decompo_dd^c_Equation_}
\eta= \6\bar\6 \alpha, \quad \alpha\in \Lambda^{p-1,q-1}(M).
\end{equation}
This statement (which is called
{global \bf $\6\bar\6$-lemma}) fails to be true
on non-K\"ahler manifolds; however, the corresponding
complex
\begin{equation}\label{_BC_intro_complex_Equation_}
\cdots\arrow\Lambda^{p-1,q-1}(M) \stackrel{\6\bar\6}\arrow \Lambda^{p,q}(M)
\stackrel{\6\oplus\bar\6}\arrow \Lambda^{p+1,q}(M)
\oplus\Lambda^{p,q+1}(M)\arrow\cdots
\end{equation}
is still elliptic. Its cohomology groups $H^{p, q}_{\6\bar\6}(M)$
are called {\bf the Bott--Chern cohomology groups of $M$} (see e.g.
\cite{_Teleman:cone_}). Explicitly, the Bott--Chern cohomology groups are:
\begin{multline*}
   H^{p, q}_{\6\bar\6}(M)=\\ = \frac{\ker
\bigg(\Lambda^{p,q}(M)\stackrel{\6}\arrow
\Lambda^{p+1,q}(M)\bigg)\bigcap \ker
\bigg(\Lambda^{p,q}(M)\stackrel{\bar\6}\arrow
\Lambda^{p,q+1}(M)\bigg)}{\im
\bigg(\Lambda^{p-1,q-1}(M)\stackrel{\6\bar\6}\arrow
\Lambda^{p,q}(M)\bigg)}.
\end{multline*}
For compact K\"ahler manifolds the Bott--Chern cohomology groups are
isomorphic with the Dolbeault ones; this isomorphism is
equivalent to the $\6\bar\6$-lemma. For non--K\"ahler manifolds, this
isomorphism does not hold anymore.

For Morse--Novikov cohomology, a similar complex can be
defined. Consider the Hodge components of the
Morse--Novikov differential $d_\theta := d-\theta$:
$d_\theta = \6_\theta+\bar\6_\theta$,
with $\6_\theta=\6-\theta^{1,0}$ and
$\bar\6_\theta=\bar\6-\theta^{0,1}$.
Locally, the Morse--Novikov complex becomes the
de Rham complex after a change $\eta\mapsto \psi\eta$,
where $\psi=e^f$, with $f$ a function which satisfies
$df=\theta$. Therefore, the  complex
\begin{equation}\label{_BC_MN_intro_complex_Equation_}
\Lambda^{p-1,q-1}(M) \stackrel{\6_\theta\bar\6_\theta}\arrow
\Lambda^{p,q}(M)
\stackrel{\6_\theta\oplus\bar\6_\theta}\arrow \Lambda^{p+1,q}(M)
\oplus\Lambda^{p,q+1}(M)
\end{equation}
is also elliptic. Its cohomology groups
$H^{p, q}_{\6_\theta\bar\6_\theta}(M)$
are called {\bf the Bott--Chern  cohomology} groups of the Morse--Novikov
complex.

It is possible to express $H^{1,1}_{\6_\theta\bar\6_\theta}(M)$
in terms of Morse--Novikov cohomology of $M$ and
holomorphic cohomology of a flat line bundle $L$,
with monodromy determined from $\theta$ (see
\ref{_BC_via_MN_Theorem_}).

The cohomology class
$[\omega]\in H^{1,1}_{\6_\theta\bar\6_\theta}(M)$
of an LCK-form $\omega$ is called {\bf the Bott--Chern  class of
$M$}.

\smallskip

To summarize: with every LCK-manifold, we associate
three cohomological invariants:
the Lee class $[\theta]\in H^1(M)$, the
Morse--Novikov class $[\omega]\in H^2_\theta(M)$,
and the Bott--Chern  class
$[\omega]\in H^{1,1}_{\6_\theta\bar\6_\theta}(M)$.
Notice that the Morse--Novikov class can be
reconstructed from the Bott--Chern  class.

In different situations, these three
cohomological invariants play the role
of the K\"ahler class. However, unlike the
K\"ahler class, the Morse--Novikov and
Bott--Chern  classes can be zero (the Morse--Novikov class  vanishes
for compact Vaisman manifolds, see Section
\ref{_MN_Section_}). On the
other hand, if $M$ is compact and non-K\"ahler,
the Lee class $[\theta]$ cannot vanish, and
if $M$ admits a K\"ahler structure, all LCK-structures
on $M$ are globally conformally K\"ahler and have $[\theta]=0$
(a result proved by Vaisman, see \cite[Theorem 2.1]{_Dragomir_Ornea_}).

It is not clear whether the global $\6\bar\6$-lemma
is true for Morse--Novikov complex, associated with the Lee
form of an LCK-structure. If it is true, then the tautological map
$H^{1,1}_{\6_\theta\bar\6_\theta}(M) \arrow H^{1,1}_\theta(M)$
is an isomorphism, and the Morse--Novikov class of an
LCK-manifold is the same as its Bott--Chern  class.

\subsection{The space of LCK-structures}

In K\"ahler geometry, the K\"ahler form $\omega$
determines the K\"ahler class $[\omega]\in H^{1,1}(M)$,
and the difference of K\"ahler forms which have the same
K\"ahler class is measured by a potential:
\[
\omega_1 - \omega = \6\bar\6\phi
\]
(this follows from the $\6\bar\6$-lemma). The space
of all K\"ahler metrics is locally modelled
on $H^{1,1}(M, \R)\times (C^\infty(M)/\const)$.
A similar local description exists for the set
of all LCK-structures on a given complex manifold,
if we fix the cohomology class $[\theta]$
of a Lee form.

Let $[\omega]\in H^{1,1}_{\6_\theta\bar\6_\theta}(M)$
be the Bott--Chern  class of an LCK-form $\omega$.
Given another LCK-form $\omega_1$, with the same
Bott--Chern  class, we can write
\begin{equation}\label{_LCK_via_pote_Equation_}
\omega_1  = \omega + d_\theta d^c_\theta \phi =\omega +
\phi (\theta \wedge I\theta + d^c \theta) - \theta \wedge d^c\phi
+ I\theta \wedge d\phi + dd^c \phi,
\end{equation}
where $d^c_\theta = - I d_\theta I = d^c - I\theta$.
Here we use implicitly the equality
\[ d_\theta d^c_\theta= -2\1\6_\theta\bar\6_\theta
\]
which is the Morse--Novikov version of
$dd^c = -2\1 \6\bar\6$.

For any real-valued function $\phi \in C^\infty (M)$,
the form \eqref{_LCK_via_pote_Equation_} satisfies
$d\omega_1 = \omega_1 \wedge \theta$, as a simple
calculation implies. If
\[ \sup_M \bigg(|dd^c\phi| + |d\phi|+ |\phi|\bigg)
\]
is sufficiently small, $\omega_1$ is also positive.
We obtained that the difference of two LCK-forms in the same Bott--Chern
class is expressed by a potential, just like in
K\"ahler case, and the set of LCK-structures is
locally parametrized by
\[
 H^{1,1}_{\6_\theta\bar\6_\theta}(M)\times
 \bigg( C^\infty(M) / \ker d_\theta d^c_\theta\bigg)
\]

With regard to the realization of cohomology classes
by LCK-forms, one could ask the questions similar
to those asked (and sometimes answered) in K\"ahler
geometry.

\hfill

\problem
Determine all 1-forms $\theta$ for which
there exists a Hermitian two-form $\omega$ having $\theta$
as its Lee form, and all the Morse--Novikov classes
which can be realized by an LCK-form.

\hfill

\problem
Let $M$ be a compact complex manifold,
admitting an LCK-metric, and $[\theta]\in H^1(M)$ its
Lee class.  Determine the set of all classes
$[\omega] \in H^{1,1}_{\6_\theta\bar\6_\theta}(M)$
such that  $[\omega]$ is the Bott--Chern class of
an LCK-structure with the Lee class $[\theta]$.

\hfill

If $\theta$ is fixed,
a sum of two LCK-forms is again LCK;
therefore, the set of possible
Morse--Novikov and
Bott--Chern  classes of LCK-structures
on a given manifold with a fixed Lee
class $[\theta]\in H^1(M)$
is a convex cone, similar to
a K\"ahler cone.

In algebraic geometry, one often finds
all kind of geometric properties of a K\"ahler
manifold encoded in the shape of its K\"ahler cone.
One would expect that the LCK-cones defined above
would be just as important.

\subsection{Potentials on coverings of LCK-manifolds}

An important special case of an LCK-manifold
is a manifold with {\it vanishing} Bott--Chern  class.
In \cite{_OV:_Potential_}, we studied the LCK-manifolds
admitting an embedding into a Hopf manifold. We have
proven a theorem, which can be stated as follows
(using the language developed in the present paper).

Recall that {\bf a linear Hopf manifold} is a complex
manifold of the form $H_A:= (\C^n\backslash 0)/\Gamma_A$,
where the group $\Gamma_A\cong \Z$ is generated by
$x\arrow A(x)$, where $A$ is a linear operator
with all eigenvalues $\{\alpha_i\}$ satisfying
$|\alpha_i| < 1$. If, in addition,
$A$ can be diagonalized, $H_A$ is called {\bf a diagonal
Hopf manifold}. It is easy to see that $H_A$
is homeomorphic to $S^{2n-1}\times S^1$. Since
this manifold has $b_1(H_A)$ odd, it cannot be
K\"ahler. However, $H_A$ is locally conformally
K\"ahler (see \cite{_Gauduchon_Ornea_}).
In \cite{_OV:_Potential_} we have studied LCK-manifolds
which admit a holomorphic (but not necessarily isometric)
embedding into a linear Hopf manifold.

The main result of
\cite{_OV:_Potential_} can now be stated in terms of
Bott--Chern cohomology:

\hfill

\theorem \label{_embed_to_Hopf_intro_Theorem_}
(\cite{_OV:_Potential_})
Let $M$ be a compact LCK-manifold, $\dim_\C M \geq 3$.
Then the following statements are equivalent.
\begin{description}
\item[(i)] $M$ admits a holomorphic embedding
into a linear Hopf manifold.
\item[(ii)] $M$ admits an LCK-structure with rational Lee class
and vanishing Bott-Chern class.
\end{description}

\hfill

In the present paper, we explore the geometry of
LCK-structures with vanishing Bott--Chern  class.
We generalize \ref{_embed_to_Hopf_intro_Theorem_}
to all manifolds with vanishing Bott--Chern  class.

\hfill

\theorem\label{_LCK_w_vanishing_BC_Theorem_}
Let $M$ be a compact LCK-manifold with vanishing
Bott--Chern  class. Then $M$ admits an LCK-structure
with vanishing Bott--Chern class and rational cohomology class of
its Lee form.

\hfill

\noindent
{\bf Proof:} See
\ref{_defo_Lee_to_rational_Corollary_}. \endproof

\hfill

{}From \ref{_LCK_w_vanishing_BC_Theorem_}
and \ref{_embed_to_Hopf_intro_Theorem_},
we obtain that an LCK-manifold of complex dimension at least $3$
admits a holomorphic embedding into a Hopf manifold
if and only if it admits an LCK-structure with
vanishing Bott--Chern  class. We conjecture
that the Bott--Chern  class of an LCK-manifold
admitting a holomorphic embedding into a Hopf manifold
always vanishes.

For an important class of LCK-manifolds
a weaker form of this conjecture can be proven.
One could define a compact Vaisman manifold
as an LCK-manifold admitting a holomorphic
flow acting by conformalities
(Subsection \ref{_Vaisman_Subsection_}).
It is known (\cite{_Gauduchon_Ornea_},
\cite{_Verbitsky:Sta_Elli_}) that all
diagonal Hopf manifolds are Vaisman.
In \cite{_OV:_Immersion_}, we have shown
that any Vaisman manifold admits a holomorphic
immersion into a diagonal Hopf manifold, and in
\cite{_OV:_Potential_} we proved that
for $\dim_\C M \geq 3$ there exists
an embedding into a diagonal Hopf manifold.
In fact, the assumption $\dim_\C M \geq 3$
is not needed, because in \cite{_Belgun_},
all 2-dimensional Vaisman manifolds were
classified, and embeddability of those
into $H_A$ can be easily checked using the same
arguments as in \cite{_OV:_Potential_}.

In the present paper, we study the
LCK-structures which can be defined on a
given Vaisman manifold. We show that
the Morse--Novikov class of any such structure
vanishes (\ref{_LCK_on_Vaisman_Theorem_}),
and the Lee class is rational.


\section{Locally conformally K\"ahler geometry}

In this section we give the necessary definitions and
properties of locally conformally K\"ahler (LCK)
manifolds.

In what follows, $M$ will denote a connected, smooth
manifold of real dimension $2n$; $I$ will be an integrable
complex structure. For a Hermitian metric $g$, we denote
with $\nabla^g$  its Levi-Civita connection and with
$\omega$ its fundamental two-form defined as
$\omega(X,Y)=g(IX,Y)$.

\subsection{LCK manifolds}

\definition
A complex manifold $(M,I)$ is LCK if it admits
a K\"ahler covering $(\tilde M, \tilde \omega)$,
such that the covering group acts by holomorphic
homotheties.

\smallskip

Equivalently, there exists on $M$ a \emph{closed} 1-form
$\theta$, called \emph{the Lee form}, such that $\omega$
satisfies the integrability condition:
$$d\omega=\theta\wedge\omega.$$

Clearly, the metric $g:=\omega(\cdot, I\cdot)$ on $M$ is
locally conformal to some K\"ahler metrics and its lift to
the K\"ahler cover in the definition is globally conformal
to $h$.

To an LCK manifold one associates the \emph{weight bundle}
$L_{\mathbb{R}}\longrightarrow M$.
It is a real line bundle associated to the
representation\footnote{ In conformal geometry, the weight
bundle usually corresponds to $\mid\det
A\mid^{\frac{1}{2n}}$. For LCK-geometry, $\mid\det A\mid^{\frac{1}{n}}$
is much more convenient.}
\[
\mathrm{GL}(2n,\RR)\ni A\mapsto \mid \det A\mid^{\frac{1}{n}}.
\]
The Lee
form induces a connection in $L_{\mathbb{R}}$ by the
formula $\nabla=d- \theta$. $\nabla$  is
associated to the Weyl covariant derivative (also denoted
$\nabla$) determined on $M$ by the
LCK metric and the Lee form (the Weyl covariant derivative
is uniquely defined by the properties $\nabla I=0$,
$\nabla g=\theta\otimes g$; in this context, $\theta$ is
called the Higgs field).
As $d\theta=0$, then $\nabla^2= d\theta=0$, and hence
$L_{\mathbb{R}}$ is flat.

The complexification of the weight bundle will be denoted by $L$. The Weyl
connection extends naturally to $L$ and its
$(0,1)$-part endows $L$ with a holomorphic structure. On the other hand,
as $L$ is flat, one can pick a nowhere degenerate section $\lambda$
satisfying
$$\nabla(\lambda)=\lambda\otimes (-\theta).$$
Accordingly, one chooses a Hermitian structure on $L$ such that
$\mid\lambda\mid =1$ and considers the associated Chern connection. Then
one proves:

\hfill

\theorem \cite[Theorem 6.7]{_Verbitsky_vanishing_}   The curvature of the
Chern
connection on $L$ with respect to the above holomorphic and Hermitian
structure is $-2\sqrt{-1}d^c\theta$.

\subsection{Vaisman manifolds}
\label{_Vaisman_Subsection_}

\definition
A Vaisman manifold is an LCK manifold with $\nabla^g$-parallel Lee form.

\hfill

On a Vaisman manifold, the Lee field $\theta^\sharp$ is real holomorphic and Killing (see
\cite{_Dragomir_Ornea_}). On compact manifolds, this statement can be
taken as a definition:

\hfill

\theorem\cite{_Kamishima_Ornea_} A compact LCK manifold admits a LCK
metric with parallel Lee form if and only if its Lie group of holomorphic
conformalities has a complex one-dimensional Lie subgroup,
acting non-isometrically on its K\"ahler covering.

\hfill

The structure of compact Vaisman manifolds is now fully understood:

\hfill

\theorem\cite{_OV:_Structure_}\label{struct} Let $(M,I,g)$ be a compact
Vaisman manifold. Then:

1) The monodromy of $L_\mathbb{R}$ is $\mathbb{Z}$.

2) $(M,g)$ is a suspension over $S^1$ with fibre a compact Sasakian
manifold $(W,g_W)$. Moreover, $M$ admits a conic K\"ahler covering
$(W\times \mathbb{R}_+, t^2g_W+dt^2)$ such that  the covering group is
an infinite cyclic group, generated by the
transformation $(w,t)\mapsto (\varphi(w), qt)$
for some Sasakian automorphism $\varphi$ and $q\in\mathbb{Z}$.

\hfill

The typical example of a compact Vaisman manifold is the
diagonal Hopf manifold
$H_A:=\C^n/ \langle A\rangle$ with
$A=\text{diag}(\al_i)$, with $|\alpha_i|<1$. The complex structure is the
projection of the standard one of $\C^n$ and the LCK
metric is constructed as follows.

Let $C>1$ be a constant. Then  one constructs on $\C^n$ the potential
$$\f(z_1,\ldots,z_n)=\sum \abs{z_i}^{\be_i}, \quad
\be_i=\log_{\abs{\al_i}^{-1}}C$$
which is acted on by  $A$ as follows: $A^*\f=C^{-1}\f$.

Hence:

$\displaystyle \Omega = \sqrt{-1} \partial\bar\partial \f$ is K\"ahler and
$\displaystyle \Gamma\cong \ZZ$ acts by holomorphic homotheties with
respect to it.

The Lee field is : $\displaystyle \theta^\sharp=-\sum z_i\log
\abs{\al_i}\partial z_i$ and one can see it is parallel.

Note than diagonal Hopf manifolds are generalizations of the rank 1 Hopf
surfaces.

Other examples (in fact, the whole list) of compact complex surfaces with
Vaisman structure are given in \cite{_Belgun_}.

Non-Vaisman LCK manifolds are some of the Inoue surfaces (cf.
\cite{_Tricerri_}, \cite{_Belgun_}) and their generalizations to higher
dimensions (\cite{_Oeljeklaus_Toma_}), the rank 0 Hopf surfaces
(\cite{_Gauduchon_Ornea_}).

\subsection{LCK manifolds with potential}

Not only the Vaisman metric of the Hopf manifold can be constructed out of
a potential on the K\"ahler covering, but the K\"ahler metric on the
universal cover of any Vaisman manifold has a global potential given by
$\mid\theta^\sharp\mid^2$.

A wider class of LCK manifolds, strictly containing the Vaisman ones, is
the following:

\hfill

\definition \label{_LCK_w_pote_Definition_}
(\cite{_OV:_Potential_})
A complex manifold $(M,I)$ is \emph{LCK  with  potential} if it admits a
K\"ahler cover $(\tilde M, \Omega)$ with global potential $\f:\; \tilde M
\rightarrow \RR_{+}$
satisfying the following conditions:

(i) $\f$ is proper (i.e. it  has compact level sets).

(ii) The monodromy map $\tau$ acts on $\f$ by multiplication with a
constant: $\tau (\f)=const \cdot \f$.

\hfill

On compact manifolds, the properness of the potential is equivalent (cf.
\cite{_OV:_Potential_}) to the deck group being isomorphic to $\Z$ (a
condition satisfied by compact Vaisman manifolds).

In \cite{_OV:_Immersion_} we showed that there exist deformations which preserve the Vaisman class on compact manifolds. Moreover, on compact manifolds, one can always deform a Vaisman structure into  a quasi--regular one.
Using a similar argument, one proves:

\hfill

\theorem (\cite{_OV:_Potential_}) The class of compact LCK manifolds with
potential is stable under small deformations.

\hfill

As a consequence, one sees that the Hopf manifold $(\CC^N\setminus
0)/\Gamma$, where now $\Gamma$ is a cyclic group generated by a
non-diagonal linear operator, is LCK with potential. This is a
generalization of the (non--Vaisman)  rank 0 Hopf surface.

The main property of LCK manifolds with potential is that they satisfy an
embedding theorem similar to the Kodaira-Nakano one in K\"ahler geometry:

\hfill

\theorem (\cite{_OV:_Potential_}) Any compact LCK manifold with potential
of complex dimension at least $3$ can be holomorphically embedded in a
Hopf manifold.  Moreover, a compact Vaisman manifold of complex dimension
at least $3$ can be holomorphically embedded in a diagonal
Hopf manifold.

\hfill


\section{Morse--Novikov cohomology}
\label{_MN_Section_}


\subsection{Morse--Novikov complex and cohomology of local systems}

Let $M$ be a smooth manifold, and $\theta$ a closed 1-form
on $M$. Denote by $d_\theta:\; \Lambda^i(M)\arrow \Lambda^{i+1}(M)$
the map $d-\theta$. Since $d\theta=0$, $d_\theta^2=0$.

Consider the complex
\[
\Lambda^0(M)\stackrel {d_\theta} \arrow \Lambda^{1}(M)
\stackrel {d_\theta} \arrow\Lambda^{2}(M) \stackrel {d_\theta} \arrow ...
\]
This complex is called {\bf  the Morse--Novikov complex},
(see e.g.  \cite{_Pajitnov_}, \cite{_Ranicki_},
\cite{_Millionschikov_})
and its cohomology {\bf the Morse--Novikov cohomology}.
In Jacobi and locally conformal symplectic  geometry, this object is called
{\bf Lichnerowicz-Jacobi, or Lichnerowicz cohomology}, motivated by
Lichnerowicz's
work \cite{_Lichnerowicz_} on Jacobi manifolds
(see e.g. \cite{_Leon_Lopez_} and \cite{Banyaga}).

For an early introduction to this subject, we refer to
\cite{_Pajitnov_}.

A closed 1-form $\theta$ defines a flat connection
$d-\theta$ on the trivial bundle $M\times C^\infty (M)$. A sheaf of
$d-\theta$-closed functions on $M$ is obviously locally
trivial, and hence it defines a local system. Its monodromy is
associated with the character $\chi:\pi_1(M) \arrow \R^{>0}$
given by the exponent $e^\theta \in H^1(M, \R^{>0})$,
considered as an element of $\R^{>0}$-valued cohomology.

The cohomology of this local system is equal to the
Morse--Novikov cohomology, as we shall see presently.

Let $L$ be a real local system on $M$, that is,
a locally trivial sheaf on $M$ locally modeled on
$M\times\R$. Assume also that $L$ is oriented, that
is, its monodromy lies in $\R^{>0}\subset \R^*$.
Then the vector bundle $L_\R = L\otimes_\R C^\infty (M)$
is trivial. The structure of a local system induces
on $L_\R$ a flat connection $\nabla$ (see \cite{_Hartshorne_}). Choose a
trivialization, that is, a nowhere degenerate section
$\xi\in L_\R$, and let
\[
\nabla_\xi:\; L_\R \arrow L_\R \otimes C^\infty( M)
\]
be the trivial connection on $L_\R$ mapping $\xi$ to 0.
Then $\theta= \nabla-\nabla_\xi$ is a 1-form
with values in $\End_{C^\infty( M)} (L_\R) = C^\infty (M)$.
Since $\nabla^2=0$, $\theta$ is closed. The following
elementary claim is well known.

\hfill

\proposition\label{loc_syst=mn}
\cite{Novikov}
The cohomology of the local system $L$
is naturally identified with the
cohomology of the Morse--Novikov
complex $(\Lambda^*(M), d_\theta)$.

\endproof

\subsection{Morse--Novikov cohomology of LCK-manifolds}

It is well--known that the LCK form represents a class
in the Morse--Novikov cohomology. Indeed, let $(M, \omega)$
be an LCK-manifold, with $\theta$ its Lee form.
Since $d\omega = \theta\wedge\omega$, we have
$d_{\theta}\omega=0$. Therefore, $\omega$
represents a cohomology class in the complex
$(\Lambda^*(M), d_{\theta})$

\hfill

\definition\label{_Morse_N_class_Definition_}
The Morse--Novikov cohomology class
$[\omega]$ of $\omega$ is called {\bf the
Morse--Novikov class of $M$}.

\hfill

This notion is similar to the notion of a K\"ahler class of
a K\"ahler manifold.

Morse--Novikov cohomology for locally conformally symplectic manifolds
was first considered in \cite{_Guedira_Lichnerowicz_} where it was
proven to vanish in the top dimension. Then Vaisman,
\cite{_Vaisman:80_}, studied it under the name of "adapted
cohomology" on LCK manifolds and
identified it with the cohomology with values in the sheaf of germs of
smooth $d_\theta$-closed functions. Later, it was proven in
\cite{_Leon_Lopez_} that it vanishes in all dimensions for compact locally
conformally symplectic manifolds which admit a compatible Riemannian metric
for which the Lee form is parallel, hence, in particular, for compact
Vaisman manifolds (see also \cite{_Ornea_}). But for
compact Vaisman manifolds, the vanishing of the
Morse--Novikov cohomology follows almost directly from the
Structure \ref{struct}.

Indeed, according to  \ref{loc_syst=mn}, this is the cohomology of the
local
system $L$. But, by  \ref{struct}, $M$
is $W\times S^1$ topologically and the monodromy of $L$ is $\ZZ$, hence
$L$ is
the pull-back $p^*L'$ of a local system $L'$ on $S^1$.

Now, the cohomology of the local system $L$ is the derived
direct image $R^iP_*(L)$,
where $P$ is a projection onto a point. By the above
remark and changing the base,
$R^iP_*(L)=R^iP_*p^*(L')=R^i({\CC}\otimes L')$, where
$\CC$ is viewed as  a trivial local system. From the
K\"unneth formula it follows that  $R^i(\CC)$ is a locally
constant sheaf on $S^1$, with fiber $H^*(W)$. By the Leray
spectral sequence of composition, the hypercohomology of
the complex of sheaves $R^*(\CC\otimes L')$ converges to
$R^iP_*(L)$. Finally, each $R^i(\CC\otimes L')$ has zero
cohomology, being a locally trivial sheaf on $S^1$
with non-trivial constant monodromy.
Therefore this spectral sequence vanishes in
$E_2$. It then converges to zero.

\hfill

\remark There exist compact LCK manifolds with non--vanishing
Morse--Novikov class.
Indeed, it is proven in \cite{Banyaga} that the compact $4$--dimensional LCK
solvmanifold constructed in \cite{Andres_Cordero_Fernandez_Mencia} has
non--vanishing Morse--Novikov class. On the other hand, it is shown in
\cite{Kamishima} that this solvmanifold is biholomorphic with an Inoue
surface which
is known, \cite{_Belgun_}, to not admit any Vaisman metric.


\section{Bott--Chern class of an LCK-form}


\subsection{LCK structures up to potential}

\definition
Let $(M, \omega_1, \theta)$ and $(M, \omega_2, \theta)$ be two
LCK-structures on the same compact manifold $M$, having the same Lee
form $\theta$. These structures are called {\bf equivalent up to a
potential} if the following conditions hold.

\begin{description}
\item[(i)]
Consider a covering $\tilde M$ where $\theta$ becomes exact, $\theta
= d f$, and let $\tilde \omega_i= e^{-f}\omega_i$ be the
corresponding K\"ahler forms on $\tilde M$. Then $\tilde \omega_1 -
\tilde \omega_2=\6\bar\6\phi$, for some smooth function $\phi:\;
\tilde M \arrow \R$.

\item[(ii)] Let $\Gamma$ be the deck transformation
(also known as monodromy) group of the covering $\tilde M \arrow M$,
and $\chi:\; \Gamma \arrow \R^{>0}$ the character corresponding to
$e^{f}$,
\[
\gamma^* e^f = \chi(\gamma) e^{f}, \ \ \ \forall \gamma\in \Gamma.
\]
Then $\phi$ has the same automorphy:
\[
\gamma^* \phi = \chi(\gamma) \phi, \ \ \ \forall \gamma\in \Gamma.
\]
A K\"ahler potential $\phi$ satisfying these
automorphy conditions is called {\bf the automorphic
potential} for the LCK metric.
\end{description}

\hfill

Given an LCK-structure on a compact complex manifold, it is very
easy to construct other LCK-structures, equivalent up to potential.
Let $(M, \omega)$ be an LCK-manifold, $\tilde\omega$ the natural
K\"ahler form on its covering $\tilde M$. Consider a function $v:\;
\tilde M \arrow \R$ which satisfies $|\nabla^2 v|_{L^\infty}<
\epsilon$, and has the same automorphy condition
\begin{equation}\label{_auto_pote_Equation_}
\gamma^* v = \chi(\gamma) v, \ \ \ \forall \gamma\in \Gamma.
\end{equation}
If $\epsilon< 1$, the form $\tilde\omega_1:= \tilde
\omega+\6\bar\6 v$ is K\"ahler, with the same automorphy properties
as $\tilde \gamma$. Therefore, $\tilde \omega_1$ induces an
LCK-structure on $\tilde M$, obviously, equivalent to
$(M,\omega,\theta)$. All LCK-structures equivalent to a given one
are obtained this way.

\hfill

\remark\label{_LCK_det_by_class_and_fu_Remark_} Denote by
$C^\infty(\tilde M)_\chi$ the space of smooth functions satisfying
\eqref{_auto_pote_Equation_}. We have shown that the set of
LCK-structures equivalent to a given one is an open convex cone in
an affine space modeled on a quotient of $C^\infty(\tilde M)_\chi$
by the kernel of $\6\bar\6$.

\subsection{Bott--Chern cohomology}

The set of equivalence classes of LCK-structures with a given
monodromy can be described in terms of cohomology.

\hfill

\definition
Let $M$ be a complex manifold, $\dim_\C M=n$, $0<p,q <n$ integer
numbers, and $L$ a complex line bundle with flat connection.
Consider the complex
\begin{equation}\label{_BC_complex_Equation_}
\arrow\Lambda^{p-1,q-1}(M,L) \stackrel{\6\bar\6}\arrow \Lambda^{p,q}(M,L)
\stackrel{\6\oplus\bar\6}\arrow \Lambda^{p+1,q}(M,L)
\oplus\Lambda^{p,q+1}(M,L)\arrow
\end{equation}
where $\6$, $\bar \6$ denote the (1,0) and (0,1)-parts of the
connection operator $\nabla:\; \Lambda^i(M, L)
\arrow\Lambda^{i+1}(M, L)$. The cohomology of
\eqref{_BC_complex_Equation_} is called {\bf the Bott--Chern
cohomology of $M$ with coefficients in $L$}, denoted by
$H^{p,q}_{\6\bar\6}(M,L)$. It is well known that
\eqref{_BC_complex_Equation_} is elliptic, hence
$H^{p,q}_{\6\bar\6}(M,L)$ is finite-dimensional.

\hfill

\definition
Let $(M,\omega,\theta)$ be an LCK-manifold, and $L$ its weight
bundle, that is, a trivial complex line bundle with the flat
connection $d-\theta$. Consider $\omega$ as a closed $L$-valued
(1,1)-form on $M$. Its Bott--Chern class $[\omega]\in
H^{1,1}_{\6\bar\6}(M,L)$ is called {\bf the Bott--Chern class of the
LCK-manifold}.

\hfill

Clearly, the vanishing of $[\omega]\in H^{1,1}_{\6\bar\6}(M,L)$
implies the existence of an automorphic potential for $M$. Hence:

\hfill

\proposition If the Bott--Chern class of an
LCK-manifold $M$ vanishes and the
monodromy of $L$ is
$\Z$,
then $M$ is LCK with potential.

\hfill

Directly from the definition we have:

\hfill

\proposition Let $M$ be a complex manifold, and $\omega_1$ and
$\omega_2$ LCK-metrics having the same Lee form $\theta$. Then the
following conditions are equivalent:
\begin{description}
\item[(i)] The Bott--Chern classes of $\omega_1$, $\omega_2$ are equal.
\item[(ii)] The LCK-structures $\omega_1$ and $\omega_2$
are equivalent up to a potential.
\end{description}

\hfill

We obtain a remarkable analogy between the K\"ahler manifolds and
LCK-manifolds. A K\"ahler structure on a complex manifold is
determined by a K\"ahler class in $H^{1,1}(M)$ and a choice of a
K\"ahler metric in this K\"ahler class. The latter is obtained by
choosing an element in a cone locally modeled on
$C^\infty(M)/\const$.

An LCK-structure on a given complex manifold with a prescribed
conformal structure is determined by a Bott--Chern class and a
choice of an LCK-metric with a prescribed Bott--Chern class.
A metric with prescribed Bott--Chern  class
metric is obtained by choosing an element in a  cone locally modeled
on $C^\infty(\tilde M)_\chi/\ker(\6\bar\6)$ (see
\ref{_LCK_det_by_class_and_fu_Remark_}).

\subsection{Bott--Chern classes and Morse--Novikov cohomology}

The holomorphic cohomology of a bundle can be realized as cohomology
of a complex
\begin{equation}\label{_Dolbeault_Equation_}
C^\infty(L) \stackrel{\bar \6} \arrow \Lambda^{0,1}(M,L)
\stackrel{\bar \6} \arrow  \Lambda^{0,2 }(M,L)  \stackrel{\bar \6}
\arrow \cdots
\end{equation}
If $L$ is equipped with a flat connection, $\6:\;
\Lambda^{0,1}(M,L)\arrow \Lambda^{1,1}(M,L)$ induces a map
\begin{equation}\label{_holo_coho_to_BC_Equation_}
H^1({\cal L}) \stackrel \6\arrow H^{1,1}_{\6\bar\6}(M,L)
\end{equation}
{}from the holomorphic cohomology of the underlying
holomorphic bundle (denoted as ${\cal L}$)
to the Bott--Chern cohomology. The complex
\begin{equation}\label{_Dolbeault_1,0_Equation_}
C^\infty(L) \stackrel{\nabla^{1,0}} \arrow \Lambda^{1,0}(M,L)
\stackrel{\nabla^{1,0}} \arrow  \Lambda^{2,0}(M,L)
\stackrel{\nabla^{1,0}}\arrow\cdots
\end{equation}
computes the holomorphic cohomology of a bundle ${\cal L}'$
with a holomorphic structure defined by the complex conjugate
of the $\nabla^{1,0}$-part of the connection. When
the bundle $L$ is real, we have $ {\cal L}\cong {\cal L}'$.
Then the cohomology of the complex
\eqref{_Dolbeault_1,0_Equation_} is naturally identified
with $\overline{H^*({\cal L})}$.  The map
$\bar\6:\;\Lambda^{1,0}(M,L)\arrow \Lambda^{1,1}(M,L)$
defines a homomorphism
\begin{equation}\label{_antihoholo_coho_to_BC_Equation_}
\overline{H^1({\cal L})} \stackrel {\bar\6}\arrow H^{1,1}_{\6\bar\6}(M,L)
\end{equation}
which is entirely similar to \eqref{_holo_coho_to_BC_Equation_}.

\hfill

The following result allows one to compute
the Bott--Chern  classes in terms of holomorphic
cohomology and Morse--Novikov cohomology.

\hfill

\theorem\label{_BC_via_MN_Theorem_}
Let $M$ be a complex manifold,
$L_\R$ a trivial real line bundle with flat connection $d-\theta$,
where $\theta$ is a real closed 1-form. Denote by $L$ its
complexification, and let $\cal L$ be the underlying holomorphic
bundle. Then there is an exact sequence
\begin{equation}\label{_MN_BC_exact_Equation_}
H^1({\cal L})\oplus \overline{H^1({\cal L})}
\stackrel {\6+\bar\6}\arrow H^{1,1}_{\6\bar\6}(M,L)
\stackrel\nu\arrow H^2_\theta(M)
\end{equation}
where $H^2_\theta(M)$ is the  Morse--Novikov cohomology, $\nu$
a tautological map, and the first arrow is obtained as a direct
sum of \eqref{_holo_coho_to_BC_Equation_} and
\eqref{_antihoholo_coho_to_BC_Equation_}.

\hfill

{\bf Proof:} If a $(1,1)$-form $\eta$ with coefficients in
$L$ is Morse--Novikov cohomologous to zero, we have
$\eta= d_\theta\alpha$, where $d_\theta=d-\theta$
is the corresponding differential. Taking
$(1,0)$ and $(0,1)$-parts, we obtain that
$\eta=\bar\6_\theta\alpha^{1,0} + \6_\theta\alpha^{0,1}$,
where $\6_\theta=\6- \theta^{1,0}$
and $\bar\6_\theta=\6- \theta^{0,1}$
are the Hodge components of $d_\theta$.
However, these operators are precisely those that
are used to define the first arrow in
\eqref{_MN_BC_exact_Equation_}. Moreover,
since $\eta$ has type $(1,1)$,
\[
\6_\theta\alpha^{1,0}=0, \text{\ \ and\ \ }\bar\6_\theta\alpha^{0,1}=0.
\]
Since $\alpha^{0,1}$ and $\alpha^{1,0}$  are closed
under the respective differentials, they represent
classes in the cohomology:
$[\alpha^{0,1}]\in H^1({\cal L})$ and
$[\alpha^{1,0}]\in\overline{H^1({\cal L})}$
Now the Bott--Chern  class of of $\eta$ is
obtained as $\6[\alpha^{0,1}]+\bar\6[\alpha^{1,0}]$,
hence the sequence \eqref{_MN_BC_exact_Equation_} is
exact. We proved \ref{_BC_via_MN_Theorem_}.
\endproof

\hfill

\proposition \label{_vanishing_of_bott-chern_}
Let $(M,\omega,\theta)$
be a compact LCK-manifold, $L$ the corresponding flat line bundle,
and $\cal L$ the underlying holomorphic bundle. Assume that $H^1(M,
{\cal L})=0$, and $H^2_\theta(M)=0$. Then
$H^{1,1}_{\6\bar\6}(M,L)=0$, and the K\"ahler form $\tilde \omega$
on the covering of $M$ admits an automorphic potential.

\hfill

{\bf Proof:} It follows immediately from \ref{_BC_via_MN_Theorem_}.
\endproof

\hfill

\remark
The vanishing of the Bott--Chern class is hard to control.
Hence the utility of \ref{_vanishing_of_bott-chern_} which reduces the
analysis to the more manageable Morse--Novikov cohomology of $M$ and
holomorphic cohomology of $L$.


\section{The Bott--Chern class and deformations}


We now want to determine the influence of the vanishing of the
Bott--Chern class of a compact LCK manifold. To this end, we need
some preliminaries about the automorphic potentials on the covering.

Let $\phi$ be an automorphic potential function on $\tilde M$.
As such, it can be thought of as a section of $L$. This means that
$\phi$ can be viewed as a function defined on $M$, which becomes important
in problems of approximations, as $M$ is compact and $\tilde M$
generally not. Moreover, as on $\tilde M$ one has $\theta=d\log
\phi$, the potential is uniquely determined, up to a constant, by
$\theta$.

Before stating the main result of this section, we prove a simple
technical result which was known to be true for Vaisman manifolds:

\hfill

\lemma \label{_LCK_through_Lee_Lemma_}
Let $M$ be a LCK manifold and $\tilde M$ a K\"ahler covering
on which the Lee form is exact.
Suppose the K\"ahler form of $\tilde M$
admits an automorphic potential. Then the
LCK-form of $M$ is conformally equivalent to
an LCK-form
\begin{equation}\label{lck_pot_formula}
\omega=-d^c\theta+\theta\wedge I(\theta),
\end{equation}
where $-\theta$ is the Lee form of $\omega$,
$d\omega = \omega\wedge\theta$.

\hfill

{\bf Proof:} If $\phi$ is the potential on the K\"ahler
covering $\tilde M$, the K\"ahler form $\tilde \omega$ on $\tilde M$
satisfies $\tilde \omega = dd^c\phi$. Then
$\omega:=\frac{dd^c\phi}{\phi}$ is an LCK-form on $M$
conformally equivalent to $\tilde\omega$, and the corresponding
Lee form can be found from $d\omega = \omega \wedge d\log\phi$,
giving $\theta= d\log\phi$.
Therefore,
\begin{equation*}
\begin{split}
\omega&=\frac{dd^c\phi}{\phi}=-\frac{d^cd\phi}{\phi}\\
&=-d^c\left(\frac{d\phi}{\phi}\right)+\frac{d\phi\wedge d^c\phi}{\phi^2}\\
&=-d^cd\log\phi+d\log\phi\wedge\frac{d^c\phi}{\phi}\\
&=-d^c\theta+\theta\wedge I(\theta).
\end{split}
\end{equation*}
\endproof

\hfill

\proposition\label{_close_Lee_then_pote_Proposition_}
Let $M$ be a compact LCK manifold such that the K\"ahler form on a
K\"ahler covering admits an automorphic
potential,
and let $\theta$ be its Lee form. Let $\theta'$ be any closed 1-form
sufficiently close to $\theta$  in
the metric
\[
|\theta-\theta'|_{L_1^{\infty}} = \sup_M |\theta-\theta'|
+ \sup_M |\nabla\theta-\nabla\theta'|,
\]
and let $\phi'$ be defined from $\theta'=d\log\phi'$.
Then $\omega':=\1\6_{\theta'}\bar\6_{\theta'}\phi'$
is also positive.

\hfill

{\bf Proof:} Note that in this statement the monodromy of the cover is
not assumed to be $\Z$.

Deforming $\theta$ to $\theta'$   changes the monodromy of the cover. Let
the new monodromy be
$\Gamma'$.
By definition, $\theta'=d\log\phi'$ is $\Gamma'$-invariant, hence
$a^*\log\phi'=\log\phi'+c$ for any $a\in \Gamma'$ (here $c$ is a real
constant).
This implies $a^*\phi'=\exp(c)\phi'$, and hence $\phi'$ is
$\Gamma'$--automorphic. As such, $\phi'$ is a trivialization of $L$.

Now $\1\6_{\theta'}\bar\6_{\theta'}\phi'>0$
is equivalent to $d_{\theta'}d_{\theta'}^c\phi'>0$. As a
trivialization, $\phi'$ is constant, hence the above inequality
is equivalent on $M$ with $dd^c\phi'>0$. Now $dd^c\phi'=\phi'\omega'$ and
{}from \ref{_LCK_through_Lee_Lemma_} this is equal to
$-d^c\theta'+\theta'\wedge I(\theta')$,
which is positive because $-d^c\theta+\theta\wedge I(\theta)>0$
and $\theta'$ is close to $\theta$. \endproof

\hfill

\corollary\label{_defo_Lee_to_rational_Corollary_}
Any compact LCK manifold with vanishing
Bott--Chern class admits an LCK metric with
potential, in the sense of \ref{_LCK_w_pote_Definition_}
(hence, if $\dim_\C M \geq 3$,
it is embeddable in a Hopf manifold).

\hfill

{\bf Proof:} The vanishing of the Bott--Chern class of $M$ assures
the existence of a potential on a K\"ahler covering $\tilde M\arrow
M$. The weight bundle $L$ is associated to the monodromy of this
covering and the monodromy can be a priori a subgroup
of $(\R^{>0}, \cdot)\cong (\R,+)$, which is not necessarily discrete.
Considering $L$ as a trivial line bundle with
connection $\nabla_{triv}-\theta$, we deform $L$ by
adding a small term to $\theta$ to obtain a
bundle $L'$ with monodromy $\Z$.

This is possible to do as follows. A local system on
$M$ is defined by a group homomorphism $H_1(M, \Z)\arrow \R$,
and its monodromy is $\Z$ if this map is rational.
However, each real homomorphism from $H_1(M, \Z)$ can be approximated
by a rational one. This allows one to deform $L$
into $L'$ with integer monodromy.

But deforming the monodromy amounts at
deforming the 1-form $\theta$ and this, as follows from
the formula $\theta=d\log\phi$, gives a
corresponding deformation of the potential $\phi$.
Hence we deform the pair
$(L,\phi)$ to a pair $(L',\phi')$ in which $\phi'$ is
automorphic function on $\tilde M$, with monodromy
determined by $L'$. By
\ref{_close_Lee_then_pote_Proposition_}, $\phi'$ is
plurisubharmonic, if $\theta'$ is sufficiently close to $\theta$.
By construction, $L'$ has monodromy $\Z$, hence $\phi'$
defines an LCK-metric with potential, in the sense of
\ref{_LCK_w_pote_Definition_}.
\endproof

\section{LCK-structures on compact Vaisman manifolds}

Given a closed
form $\eta$ on a manifold $M$
with an action of a connected compact group $G$,
we can average this form over the group action, obtaining
a closed form $\eta_G$. If $\eta$ is exact, $\eta_G$ is
also exact. Since $G$ is connected, it acts trivially on cohomology.
Therefore, the form $\eta$ is cohomology equivalent to $\eta_G$:
\begin{equation}\label{_coho_class_equi_Equation_}
[\eta]=[\eta_G]
\end{equation}

The same is true for $d_\theta$-closed forms,
if the form $\theta$ is $G$-invariant. Indeed,
the action of $G$ maps $d_\theta$-closed and $d_\theta$-exact
forms to $d_\theta$-closed and $d_\theta$-exact forms,
and acts trivially on the cohomology of the local
system defined by $\theta$. Therefore,
averaging a $d_\theta$-closed form $\eta$, we obtain
a $d_\theta$-closed form $\eta_G$, which is cohomology
equivalent to $\eta$.

If $\eta$ is symplectic, the form $\eta_G$
is not necessarily symplectic (this form can be degenerate).
However, if $(M, \eta)$ is a Hermitian manifold, and $G$ preserves
the complex structure on $M$, the form $\eta_G$ is positive
definite, because it is an average of positive definite forms
over a compact group. In particular, if $(M, \eta)$
is K\"ahler, then the form $\eta_G$ is a K\"ahler form too.

By this argument, it follows that the
average of an LCK-form $\omega$ over the action
of a compact group $G$ is again an LCK-form,
assuming that the Lee form of $\omega$ is
$G$-invariant.

\hfill

\theorem\label{_LCK_on_Vaisman_Theorem_}
Let $(M,J)$ be a compact complex manifold endowed with a Vaisman structure with two--form $\omega$ and Lee form $\theta$, $\dim_\C M \geq 3$. Let
$\omega_1$ be another LCK-form (not necessarily Vaisman) on $(M,J)$,
and $\theta_1$ its Lee form. Then $\theta_1$
is cohomologous with the Lee form of a
Vaisman metric, and the Morse--Novikov class
of $\omega_1$ vanishes.

\hfill

{\bf Proof:}
Denote by $\rho$ the
Lee flow corresponding to the Vaisman structure $\omega$
on $M$. It is well known that $\rho$ can be chosen
with compact leaves, giving an action
$\rho:\; S^1 \arrow \Aut(M)$ (see \cite{_OV:_Immersion_}).
Let $\theta_1$ be the Lee form of $\omega_1$.

For any closed form $\theta_1'$
in the same cohomology class as $\theta$, we can
find an LCK-form conformally equivalent
to $\omega_1$, with the Lee form equal to $\theta_1'$.
Averaging $\theta_1$ over $\rho$, we find a
1-form which is $\rho$-invariant and has the
same cohomology class, by \eqref{_coho_class_equi_Equation_}.
Hence, replacing $\omega_1$ by
a conformally equivalent form, we may assume from the beginning that
$\theta_1$ is $\rho$-invariant.

For any $t\in S^1$, $\rho(t)^*\omega_1$ satisfies
\begin{equation}\label{_rho^*_theta_inva_Equation_}
  d(\rho(t)^*\omega_1) = \rho(t)^*\omega_1\wedge \rho(t)^*\theta_1
  = \rho(t)^*\omega_1\wedge \theta_1.
\end{equation}

Averaging $\omega_1$ over $S^1$ and applying
\eqref{_rho^*_theta_inva_Equation_}, we find an
$S^1$-invariant Hermitian form $\omega_1'$ which satisfies
\[
d\omega_1' = \omega_1'\wedge \theta_1.
\]
Since the Morse--Novikov class takes values in a
cohomology group
which is $S^1$-invariant, it does not change under averaging,
and $\omega_1'$ has the same Morse--Novikov class as $\omega_1$.
Therefore, we may also assume that $\omega_1$ is $\rho$-invariant.
Denote by $G_0$ the closure of the group of holomorphic and
conformal  automorphisms of $M$ generated by $I(\theta^\sharp)$.
$G_0$ is compact, as a closed subgroup in the compact group of
conformalities of $(M, [g])$ (see also \cite{_Kamishima_Ornea_}) and
commutative because the group generated by $I(\theta^\sharp)$ is
such and taking limits preserves commutativity. Repeating the same
procedure as above, we may assume that $\theta_1$ and $\omega_1$ are
$G_0$-invariant.

Let now $\tilde M\stackrel\pi\longrightarrow M$ be a connected covering
of $M$ on which
the pullback of $\theta_1$ is exact, hence $\tilde M$ is globally
conformal K\"ahler. Let also
$\tilde \rho:\; \R \arrow \Aut(\tilde M)$
be the holomorphic flow obtained
by lifting $\rho$ to $\tilde M$.

Denote by $\tilde \omega_1$ a K\"ahler form on $\tilde M$
globally conformal to the lift of $\omega_1$.
For all $t\in \R$, the form $\tilde \rho(t)^*\tilde \omega_1$
is a K\"ahler form, conformally equivalent to $\omega_1$.
Since $\dim_\C M >2$ and $\tilde M$ is connected, the conformal factor
$\chi$ is a constant (indeed, in general, if $\alpha$, $\alpha'$ are
closed
conformal two-forms, $\alpha$ non-degenerate, and $\alpha'=f\alpha$, then
$df\wedge\alpha=0$ implies $df=0$).
In \cite{_Kamishima_Ornea_} it was shown that
if $\theta^\sharp$  and $I(\theta^\sharp)$ act conformally and
holomorphically on an LCK-manifold, and $\theta^\sharp$ cannot be
lifted to an isometry of $\tilde M$, then $M$ is Vaisman. Unless the
conformal factor $\chi$ is $1$, we may apply this result
and obtain that
$(M, \omega_1)$ is Vaisman.
It remains to prove
\ref{_LCK_on_Vaisman_Theorem_} assuming that
$\tilde \omega_1$ is $\tilde \rho$-invariant.

Let $\phi:\; \tilde M \arrow \R$ be a function defined by
$\tilde\omega_1 =  \phi\pi^* \omega_1$. Clearly,
$\pi^*\theta_1 = d\log \phi$. Since $\phi$ is
$\tilde \rho$-invariant, we have
\begin{equation}\label{_Lie_v_phi_Equation_}
\Lie_{\tilde v} \log \phi = \langle d\log \phi, {\tilde v}\rangle =0,
\end{equation}
where $\tilde v$ is the vector field generating the flow $\tilde \rho$.

>From \eqref{_Lie_v_phi_Equation_} it follows that
$\langle \theta_1, v\rangle=0$, where $v$ is the Lee field
of the Vaisman structure $\omega$. As $\theta_1$ is invariant, this means that
$\theta_1$ is a {\it basic form} associated with
the foliation $\rho$ (see \cite{_Tondeur_} for a definition
and fundamental properties of basic forms with respect to foliations).
Now, a basic 1-form on a Vaisman manifold is cohomologous
to a sum of a holomorphic and antiholomorphic
basic forms: this follows by applying the Hodge decomposition to basic
forms with respect to a transversally K\"ahler foliation.
Then $d^c\theta_1=0$.
This implies
\begin{equation}\label{_dd^c_omega_holo_Equation_}
d d^c \omega_1^{n-1} = (n-1)^2\theta_1 \wedge I(\theta_1)\wedge
\omega^{n-1}_1
\end{equation}
where $n=\dim_\C M$.
The form \eqref{_dd^c_omega_holo_Equation_} is manifestly
positive, and strictly positive unless $\theta_1$ is identically zero.
Since
\[
0 = \int_M d d^c \omega^{n-1}_1 =
\int_M (n-1)^2\theta_1 \wedge I(\theta_1)\wedge \omega^{n-1}_1,
\]
the form $\theta_1 \wedge I(\theta_1)\wedge \omega^{n-1}_1$
vanishes everywhere, hence $\theta_1$
is zero. Then $(M,\omega_1)$ is K\"ahler. This is a contradiction, since
a compact K\"ahler manifold cannot support a Vaisman structure
as compact Vaisman manifolds have odd first Betti number
(\cite{_Kashiwada_Sato_}, \cite{_Dragomir_Ornea_}). We proved
\ref{_LCK_on_Vaisman_Theorem_}. \endproof

\hfill

\remark
Notice that \ref{_LCK_on_Vaisman_Theorem_}
implies, in particular, that the weight bundle of
any LCK-structure $\omega_1$ on a Vaisman manifold has monodromy $\Z$.
Indeed, $M$ has the same monodromy as a Vaisman
manifold, and for Vaisman manifold the
weight bundle has monodromy  $\Z$, as follows from
\cite{_OV:_Structure_}.

\hfill

It is interesting to determine all
Bott--Chern  classes realized by an LCK-form on a Vaisman manifold.
In this regard, we state:

\hfill

\conjecture
Let $M$ be a Vaisman manifold, equipped with an
additional LCK-form $\omega_1$ (not necessarily Vaisman).
Then the Bott--Chern  class of $\omega_1$
vanishes; equivalently, $\omega_1$ is an
LCK-structure with potential.

\hfill

\noindent{\bf Acknowledgements:} Liviu Ornea thanks the University
of Glasgow, the Independent University and the Steklov Institute in
Moscow for hospitality during part of the work at this paper.

\hfill

{\small

\noindent {\sc Liviu Ornea\\
University of Bucharest, Faculty of Mathematics, \\14
Academiei str., 70109 Bucharest, Romania.}\\
{\sc Institute of Mathematics "Simion Stoilow" of the Romanian
Academy,\\
21, Calea Grivitei Str.
010702-Bucharest, Romania\\
\tt lornea@gta.math.unibuc.ro, \ \  Liviu.Ornea@imar.ro

\hfill

\noindent {\sc Misha Verbitsky\\
 Institute of Theoretical and
Experimental Physics \\
B. Cheremushkinskaya, 25, Moscow, 117259, Russia }\\
\tt verbit@maths.gla.ac.uk, \ \  verbit@mccme.ru
}

}

\end{document}